\documentclass[a4paper]{article}

\usepackage[margin=1in]{geometry} 

\usepackage[T1]{fontenc}
\newcommand{\changefont}[3]{
\fontfamily{#1} \fontseries{#2} \fontshape{#3} \selectfont}

\changefont{ptm}{m}{n}

\usepackage{setspace} 
\usepackage{graphicx}

\usepackage{amsfonts}
\usepackage{amsmath}
\usepackage{amssymb}
\usepackage{graphics}
\usepackage{mathrsfs}
\usepackage{color}
\newtheorem{remark}{Remark}[section]

\newtheorem{theorem}{Theorem}[section]

\newtheorem{lemma}{Lemma}[section]

\newtheorem{definition}{Definition}[section]

\long\def\symbolfootnote[#1]#2{\begingroup%
\def\thefootnote{\fnsymbol{footnote}}\footnote[#1]{#2}\endgroup} 

\begin{document}

\begin{center}
\Large \textbf{Asymptotically Unpredictable Solutions of Quasilinear Impulsive Systems with Regular Discontinuity Moments}
\end{center}

\begin{center}
\normalsize \textbf{Mehmet Onur Fen$^{1,} \symbolfootnote[1]{Corresponding Author. E-mail: monur.fen@gmail.com, onur.fen@tedu.edu.tr}$, Fatma Tokmak Fen$^2$} \\
\vspace{0.2cm}
\textit{\textbf{$^1$Department of Mathematics, TED University, 06420 Ankara, Turkey}} \\

\vspace{0.1cm}
\textit{\textbf{$^2$Department of Mathematics, Gazi University, 06560 Ankara, Turkey}} \\
\vspace{0.1cm}
\end{center}

\vspace{0.3cm}

\begin{center}
\textbf{Abstract} 
\end{center}

\noindent The notion of asymptotic unpredictability was recently introduced in (Commun. Nonlinear Sci. Numer. Simul. 134, 108029, 2024) for semiflows. Likewise unpredictable trajectories, asymptotically unpredictable ones are also capable of producing sensitivity in a dynamics, which is an indispensable feature of chaos. In the present study, we newly propose piecewise continuous asymptotically unpredictable functions, and investigate the existence and uniqueness of such solutions in a quasilinear impulsive system of differential equations comprising a term which is periodic in the time argument. The class of functions and the impulsive system under discussion admit regular discontinuity moments. Some techniques for obtaining discontinuous asymptotically unpredictable functions are additionally provided. Even though piecewise continuous unpredictable functions are asymptotically unpredictable, it is demonstrated that the converse is not true. In other words, the set of discontinuous unpredictable functions is properly contained in the set of asymptotically unpredictable ones. Appropriate examples are given with regard to all theoretical results.
\vspace{-0.2cm}

\noindent\ignorespaces

\vspace{0.3cm}
 
\noindent\ignorespaces \textbf{Keywords:} Impulsive differential equations, piecewise continuous asymptotically unpredictable functions, asymptotically unpredictable sequences, regular impulse moments.

\vspace{0.2cm}

\noindent\ignorespaces \textbf{Mathematics Subject Classification:} 34A37, 34C60

\vspace{0.6cm}

%%%%%%%%%%%%%%%%%%%%%%%%%%%%%%%%%%%%%%%%%%%%%%%%%%%%%%%%%%%%%%%%%%%%%%%%%%%%%%%%%%%%%%%%%%%%%% 

\section{Introduction} \label{sec1}

An answer to the question "Can a single trajectory generate chaos in a semiflow?" was proposed in paper \cite{AkhmetFen2016} by introducing the concept of unpredictability. The chaos type produced by an unpredictable trajectory is called Poincar\'e chaos, and its ingredients are sensitivity, transitivity and the existence of an uncountable set of motions everywhere dense and positively  Poisson stable \cite{AkhmetFen2016}. A general version of unpredictable points was discussed by Miller \cite{Miller19} for semiflows with arbitrary acting abelian topological monoids.  Topologically unpredictable points, on the other hand, were studied by Mahajan et al. \cite{Mahajan2024} for semiflows on topological spaces. Besides, the notion of Poincar\'e chaos on the arbitrary product of semiflows was investigated in \cite{Thakur21}, and developments were performed in \cite{Fen2017,Fen2018} for unpredictable solutions of differential and difference equations. Recently, it was demonstrated by Fen and Tokmak Fen \cite{Fen2024} that sensitivity in a semiflow can be attained under a weaker hypothesis compared to unpredictability. This was performed by introducing the new type of trajectory called asymptotically unpredictable. Continuous asymptotically unpredictable functions were proposed in \cite{Fen2024} taking advantage of the Bebutov dynamical system \cite{Sell1971}.

In this study we mainly focus on asymptotically unpredictable oscillations generated by systems of impulsive differential equations. In general, this type of systems are suitable for modeling processes in which abrupt changes take place \cite{Akhmet2010}, and they can be regarded as hybrid systems since they comprise both continuous and discrete dynamics \cite{Wang2020}. Impulsive models have applications in various fields such as control theory, neural networks, robotics, predator-prey systems, and epidemic dynamics \cite{Yang99}-\cite{Mondal2024}. 

We rigorously prove the existence and uniqueness of asymptotically unpredictable solutions in quasilinear systems with regular moments of impulses. This is achieved by making benefit of an impulsive system possessing an unpredictable solution, which was discussed in \cite{Fen2023}. A Gronwall type inequality for piecewise continuous functions \cite{Bainov1992} and integro-sum equations satisfied by bounded solutions are utilized. The asymptotic stability of the solutions is also taken into consideration. To investigate the model under discussion from the qualitative point of view, we provide a novel definition of asymptotically unpredictable functions with regular discontinuity moments. Such functions can be decomposed as the sum of an unpredictable function and another one converging to zero, both of which comprise discontinuities of the same type. Techniques for generating new discontinuous asymptotically unpredictable functions from a given one are also provided. It is worth noting that even though a discontinuous unpredictable function is asymptotically unpredictable, the converse is not true in general. We demonstrate the existence of discontinuous asymptotically unpredictable functions which are not unpredictable. In other words, the set of discontinuous asymptotically unpredictable functions properly includes the set of unpredictable ones.

Unpredictable solutions of impulsive systems were discussed in the papers \cite{Fen2023, Akhmet2022}. Our results are different from these studies since we take into account asymptotically unpredictable solutions instead of unpredictable ones. Owing to the existence of non-unpredictable piecewise continuous asymptotically unpredictable functions, the investigation of asymptotic unpredictability in impulsive systems cannot be reduced to the case of unpredictable solutions. This approves the novelty of the present article. On the other hand, a method for the construction of discontinuous unpredictable functions was given in paper \cite{Zhamanshin2023}. Differently from \cite{Zhamanshin2023}, we focus on the construction of asymptotically unpredictable functions.

The remainder of this research is organized in the following way. 
In Section \ref{sec2}, the notion of discontinuous asymptotically unpredictable functions is introduced and methods for generating such functions are provided. Moreover, the presence of non-unpredictable piecewise continuous asymptotically unpredictable functions is shown. In Section \ref{sec3}, we demonstrate the existence and uniqueness of asymptotically unpredictable solutions in quasilinear systems with regular moments of impulses. Examples of discontinuous asymptotically unpredictable functions and an impulsive system possessing such a solution are provided in Section \ref{sec4}. In particular, we construct a discontinuous asymptotically unpredictable function which is not unpredictable. Finally, Section \ref{concsec} is devoted to concluding remarks.

%%%%%%%%

\section{Piecewise Continuous Asymptotically Unpredictable Functions}\label{sec2}

Throughout the paper, $\mathbb R$, $\mathbb Z$, and $\mathbb N$ stand for the sets of real numbers, integers, and natural numbers, respectively.

If $\theta=\{\theta_k\}_{k\in\mathbb Z}$ is a strictly increasing sequence of real numbers satisfying $\left|\theta_k \right| \to \infty$ as $\left|k \right|\to\infty$, then we say that a piecewise continuous and bounded function $\varphi$ defined on $\mathbb R$ with values in $\mathbb R^m$ belongs to the set $\mathcal{BPC}(\theta)$ provided that for each integer $k$ the function $\varphi$ is continuous on the interval $(\theta_{k-1},\theta_k)$, the one sided limit $\varphi(\theta_k+)=\displaystyle \lim_{t \to \theta_k^+} \varphi(t)$ exists, and $\varphi(\theta_k)=\displaystyle \lim_{t \to \theta_k^-} \varphi(t)$ \cite{Akhmet2010}. Moreover, we will denote by $\mathcal{S}(p,\omega)$ the set of all strictly increasing sequences of real numbers $\theta=\{\theta_k\}_{k\in\mathbb Z}$ such that $\theta_{k+p} =\theta_k+\omega$ for each $k\in\mathbb Z$, where $p$ is a natural number and $\omega$ is a positive number. 

In the sequel, we utilize the Euclidean norm for vectors and the spectral norm for square matrices. The description of an unpredictable function with regular discontinuity moments is as follows.

\begin{definition} (\cite{Fen2023}) \label{impuldefn1}
Suppose that $\theta=\{\theta_k\}_{k\in\mathbb Z}$ is an element of $\mathcal{S}(p,\omega)$ for some $p\in\mathbb N$ and $\omega >0$. A function $\psi \in \mathcal{BPC}(\theta)$ is called unpredictable if there exist positive numbers $\varepsilon_0$ (the unpredictability constant), $\nu$ and sequences $\{\mu_n\}_{n\in\mathbb N}$, $\{\tau_n\}_{n\in\mathbb N}$ of real numbers both of which diverge to infinity such that 
\begin{itemize}
	\item[\textrm{i}.] for every positive number $\varepsilon$ there exists a positive number $\delta$ such that $\left\|\psi(s_1) -\psi(s_2) \right\|<\varepsilon$ whenever the points $s_1$ and $s_2$ belong to the same interval of continuity and $\left|s_1-s_2 \right|<\delta$;
	\item[\textrm{ii}.] $\left\| \psi(t+\mu_n) -\psi(t)\right\| \to 0$ as $n \to \infty$ uniformly on compact subsets of $\mathbb R$;
	\item[\textrm{iii}.] $\left\| \psi(t + \mu_n) - \psi(t) \right\| \geq \varepsilon_0$ for each $t \in \left[\tau_n-\nu, \tau_n+\nu \right] $ and $n\in\mathbb N$. 
\end{itemize}
\end{definition}
	
 We define an asymptotically unpredictable function belonging to $\mathcal{BPC}(\theta)$, where $\theta\in \mathcal{S}(p,\omega)$, in the following way.
 
\begin{definition} \label{impuldefn2}
Suppose that $\theta=\{\theta_k\}_{k\in\mathbb Z}$ is an element of $\mathcal{S}(p,\omega)$ for some $p\in\mathbb N$ and $\omega >0$. A function $\phi \in \mathcal{BPC}(\theta)$ is called asymptotically unpredictable if there exist an unpredictable function $\psi \in \mathcal{BPC}(\theta)$ and a function $\xi \in \mathcal{BPC}(\theta)$ satisfying $ \left\| \xi(t) \right\|  \to 0$ as $t \to\infty$ such that $\phi(t) = \psi(t) + \xi(t)$ for every $t\in\mathbb R$.
\end{definition}

The subsequent assertion provides a technique to obtain new asymptotically unpredictable functions with discontinuities from a given one.

\begin{lemma} \label{disconlemma1}
Let $\phi\in \mathcal{BPC}(\theta)$ be an asymptotically unpredictable function, where $\theta=\{\theta_k\}_{k\in\mathbb Z}$ is a sequence in $\mathcal{S}(p,\omega)$ for some $p\in\mathbb N$ and $\omega >0$. Then, for every $c \in\mathbb R^m$ and for every non-singular matrix $\Omega \in \mathbb R^{m\times m}$ the function $\widetilde{\phi} \in \mathcal{BPC}(\theta)$ defined by $\widetilde{\phi}(t) = \Omega\phi(t)+c$ is also asymptotically unpredictable.
\end{lemma}

\noindent \textbf{Proof.} Owing to the asymptotic unpredictability of $\phi$, the equation $\phi(t) = \psi(t) + \xi(t)$ is fulfilled for $t\in\mathbb R$, where $\psi$ and $\xi$ are functions in $\mathcal{BPC}(\theta)$ such that the former is unpredictable whereas the latter has limit zero as $t\to \infty$. In this case, the items (i)-(iii) of Definition \ref{impuldefn1} are fulfilled for the function $\psi$ for some positive numbers $\varepsilon_0$, $\nu$ and sequences $\{\mu_n\}_{n\in\mathbb N}$, $\{\tau_n\}_{n\in\mathbb N}$ of real numbers both of which diverge to infinity. Let us fix $c\in\mathbb R^m$ and a non-singular matrix $\Omega\in \mathbb R^{m\times m}$. Define the functions $\widetilde{\psi}$ and $\widetilde{\xi}$ in $\mathcal{BPC}(\theta)$ respectively by means of the equations $\widetilde{\psi}(t)= \Omega\psi(t)+c$ and $\widetilde{\xi}(t)=\Omega \xi(t)$, $t\in\mathbb R$. We have $\left\|\widetilde{\xi}(t) \right\| \to 0$ as $t\to\infty$ since $\left\| \widetilde{\xi}(t)\right\| \leq \left\| \Omega\right\|\left\| \xi(t)\right\|$.

Now, let a positive number $\varepsilon$ be given. There is a positive number $\delta$ such that if the points $s_1$ and $s_2$ belong to the same interval of continuity and $\left|s_1-s_2 \right|<\delta$, then $$\left\|\psi(s_1) -\psi(s_2) \right\|<\dfrac{\varepsilon}{\left\| \Omega\right\|}.$$ Therefore, 
$$
\Big\|\widetilde{\psi}(s_1) - \widetilde{\psi}(s_2) \Big\| \leq \left\|\Omega \right\| \left\| \psi(s_1) -\psi(s_2) \right\|  < \varepsilon, 
$$ whenever $s_1$ and $s_2$ belong to the same interval of continuity and $\left|s_1-s_2 \right|<\delta$.

Next, suppose that $\mathcal{C}$ is a compact subset of the real axis. Since  $\left\| \psi(t+\mu_n) -\psi(t)\right\| \to 0$ as $n \to \infty$ uniformly on $\mathcal{C}$, utilizing the inequality
$$\left\| \widetilde{\psi}(t+\mu_n) - \widetilde{\psi}(t)\right\|\leq \left\| \Omega\right\|  \left\| \psi(t+\mu_n) -\psi(t)\right\|$$ one can confirm that $\left\| \widetilde{\psi}(t+\mu_n) - \widetilde{\psi}(t)\right\| \to 0$ as $n \to \infty$ uniformly on $\mathcal{C}$. On the other hand, for each $t \in \left[\tau_n-\nu, \tau_n+\nu \right] $ and $n\in\mathbb N$ we have
$$\left\| \widetilde{\psi}(t + \mu_n) - \widetilde{\psi}(t) \right\| \geq \displaystyle \frac{1}{\left\|\Omega^{-1} \right\|} \left\| \psi(t + \mu_n) - \psi(t) \right\| \geq \frac{\varepsilon_0}{\left\|\Omega^{-1} \right\|}.$$
Thus, the function $\widetilde{\psi}$ is unpredictable with unpredictability constant $\varepsilon_0/\left\|\Omega^{-1}\right\| $. 

Making use of the equation $\widetilde{\phi}(t)= \widetilde{\psi}(t) + \widetilde{\xi}(t)$, $t\in\mathbb R$, one can deduce that $\widetilde{\phi}$ is asymptotically unpredictable. $\square$ 

Another result for piecewise continuous asymptotically unpredictable functions is as follows.

\begin{lemma} \label{disconlemma2}
Let $\theta=\{\theta_k\}_{k\in\mathbb Z}$ be an element of $\mathcal{S}(p,\omega)$ for some $p\in\mathbb N$ and $\omega >0$, and suppose that $\phi \in \mathcal{BPC}(\theta)$ is an asymptotically unpredictable function. If $\eta \in \mathcal{BPC}(\theta)$ is a function such that the limit $\displaystyle \lim_{t\to\infty} \eta(t)$ exists, then the function $\widetilde{\phi} \in \mathcal{BPC}(\theta)$ defined by $\widetilde{\phi}(t) = \phi(t) + \eta(t)$ is also asymptotically unpredictable.
\end{lemma}

\noindent \textbf{Proof.} According to Definition \ref{impuldefn2}, one can decompose $\phi$ as the sum of two functions $\psi$ and $\xi$ in $\mathcal{BPC}(\theta)$ such that $\psi$ is unpredictable and $\displaystyle \lim_{t \to\infty} \left\| \xi(t)\right\|=0$. Suppose that $\displaystyle \lim_{t\to\infty} \eta(t)=\eta_0$ for some $\eta_0 \in \mathbb R^m$, and denote $\widetilde{\xi}(t) = \xi(t) + \eta(t) -\eta_0$. The inequality $$\left\| \widetilde{\xi}(t)  \right\| \leq \left\| \xi(t)\right\| + \left\| \eta(t) -\eta_0\right\|, \, t \in\mathbb R$$ implies that  $\left\| \widetilde{\xi}(t) \right\| \to 0$ as $t\to\infty$. Hence,  the function $\overline{\phi} \in \mathcal{BPC}(\theta)$ satisfying $\overline{\phi}(t) = \psi(t) + \widetilde{\xi}(t)$ is asymptotically unpredictable. Because the equation
$$\widetilde{\phi}(t) = \overline{\phi}(t) + \eta_0$$ holds, $\widetilde{\phi}$ is asymptotically unpredictable by Lemma  \ref{disconlemma1}. $\square$

The next lemma indicates that the asymptotic unpredictability feature of a piecewise continuous function is permanent under a shifting in the time argument.

\begin{lemma} \label{disconlemma3}
Suppose that $\phi\in\mathcal{BPC}(\theta)$ is an asymptotically unpredictable function, where the sequence $\theta=\{\theta_k\}_{k\in\mathbb Z}$ belongs to $\mathcal{S}(p,\omega)$ for some $p\in\mathbb N$ and $\omega >0$. Then, for every real number $c$, the function $\widetilde{\phi}: \mathbb R \to \mathbb R^m$ defined by $\widetilde{\phi}(t) = \phi(t+c)$ is also asymptotically unpredictable. 
\end{lemma}

\noindent \textbf{Proof.} Let $\psi \in \mathcal{BPC}(\theta)$ be an unpredictable function and $\xi \in \mathcal{BPC}(\theta)$ be a function with $\left\| \xi(t)\right\|  \to 0$ as $t \to \infty$ such that $\phi(t) = \psi(t) + \xi(t)$, $t \in\mathbb R$. Fix an arbitrary real number $c$, and define the functions $\widetilde{\psi}:\mathbb R \to \mathbb R^m$,  $\widetilde{\xi}:\mathbb R \to \mathbb R^m$ respectively by the equations $\widetilde{\psi}(t) = \psi(t+c)$ and $\widetilde{\xi}(t)=\xi(t+c)$. The functions $\widetilde{\psi}$ and $\widetilde{\xi}$ are elements of $\mathcal{BPC} \Big( \widetilde{\theta} \,\Big) $, where $\widetilde{\theta}=\left\lbrace \widetilde{\theta}_k \right\rbrace _{k\in\mathbb Z}$ is the sequence in $\mathcal{S}(p,\omega)$ satisfying $\widetilde{\theta}_k = \theta_k-c$ for all $k\in\mathbb Z$. One can confirm that $\left\| \widetilde{\xi}(t)\right\|  \to 0$ as $t \to \infty$.

Now, suppose that $\varepsilon_0$, $\nu$ are positive numbers and $\{\mu_n\}_{n\in\mathbb N}$, $\{\tau_n\}_{n\in\mathbb N}$ are sequences of real numbers both of which diverge to infinity such that the items $(i)-(iii)$ of Definition \ref{impuldefn1} hold for the function $\psi$.

For a fixed $\varepsilon>0$, one can find a number $\delta >0$ such that for each elements $s_1,s_2$ of any one of the intervals $(\theta_{k-1},\theta_k)$, $k\in\mathbb Z$, with $|s_1-s_2|<\delta$, we have $\left\| \psi(s_1)-\psi(s_2) \right\|<\varepsilon$. If $s_1,s_2$ belong to one of the intervals $\left( \widetilde{\theta}_{k-1},\widetilde{\theta}_k \right)$, $k \in \mathbb Z$, and if $|s_1-s_2|<\delta$, then $s_1+c,s_2+c$ are elements of $\left( \theta_{k-1},\theta_k \right)$ so that
$$\left\| \widetilde{\psi}(s_1) - \widetilde{\psi}(s_2) \right\| = \left\| \psi(s_1+c) - \psi(s_2+c) \right\|<\varepsilon.$$ For that reason, item (i) of Definition \ref{impuldefn1} is valid for $\widetilde{\psi}$.

Next, we take a compact subset $\mathcal{C}$ of $\mathbb R$, and suppose that $\mathcal{C}\subseteq [a,b]$ for some real numbers $a$ and $b$ with $b>a$. Let a positive number $\varepsilon$ be given. Due to the uniform convergence of $\left\| \psi(t+\mu_n) -\psi(t)\right\|$ on the interval $[a+c,b+c]$ to $0$, there is a natural number $n_0$ such that if $n \geq n_0$, then $\left\| \psi(t+\mu_n) -\psi(t)\right\| < \varepsilon$ for every $t \in  [a+c,b+c]$. Accordingly, we have for every $t \in [a,b]$ that $$\left\| \widetilde{\psi}(t+\mu_n) -\widetilde{\psi}(t)\right\| = \left\| \psi(t+c+\mu_n) - \psi(t+c)\right\| < \varepsilon$$ whenever $n \geq n_0$. Therefore, 
$\left\| \widetilde{\psi}(t+\mu_n) - \widetilde{\psi}(t)\right\| \to 0$ as $n \to \infty$ uniformly on $\mathcal{C}$. This discussion approves item (ii) of Definition \ref{impuldefn1}.

Finally, to verify the last item of Definition \ref{impuldefn1} for $\widetilde{\psi}$, we define $\widetilde{\tau}_n=\tau_n-c$ for each $n \in\mathbb N$. The sequence $\left\lbrace \widetilde{\tau}_n\right\rbrace_{n\in\mathbb N}$ diverges to infinity since the same is true for $\{\tau_n\}_{n\in\mathbb N}$. Because  $\left\| \psi(t + \mu_n) - \psi(t) \right\| \geq \varepsilon_0$ for each $t \in \left[\tau_n-\nu, \tau_n+\nu \right] $ and $n\in\mathbb N$, one can attain that
$$\left\| \widetilde{\psi}(t + \mu_n) - \widetilde{\psi}(t) \right\| = \left\| \psi(t +c+ \mu_n) - \psi(t+c) \right\| \geq \varepsilon_0, \ t \in [\widetilde{\tau}_n -\nu, \widetilde{\tau}_n+\nu], \ n\in\mathbb N.$$
Consequently, $\widetilde{\phi} \in \mathcal{BPC} \Big( \widetilde{\theta} \,\Big)$ is asymptotically unpredictable. $\square$

One of the crucial results of the present study is given in the next assertion. It reveals that there exist discontinuous asymptotically unpredictable functions which are not unpredictable. In particular, it provides a technique for constructing such functions. For a piecewise continuous function $\varphi:\mathbb R \to \mathbb R^m$ the notation $\varphi(s+)$ stands for the right limit of $\varphi$ at $s$.

\begin{lemma} \label{disconlemma4}
Suppose that $\psi \in \mathcal{BPC}(\theta)$ is an unpredictable function such that $\displaystyle \sup_{t \in \mathbb R} \left\| \psi(t) \right\| \leq H$ for some $H >0$, where $\theta=\{\theta_k\}_{k\in\mathbb Z}$ belongs to $\mathcal{S}(p,\omega)$ for some $p\in\mathbb N$ and $\omega >0$. If $\xi \in \mathcal{BPC}(\theta)$ is a function with $\left\|\xi(t) \right\| \to 0$ as $t \to \infty$ such that either $\left\| \xi(t_0)\right\| \geq 4H$ for some $t_0 \in\mathbb R$ or $\left\| \xi(\theta_{k_0}+)\right\| > 4H$ for some $k_0 \in \mathbb Z$, then the asymptotically unpredictable function $\phi \in \mathcal{BPC}(\theta) $ satisfying $\phi(t) = \psi(t) + \xi(t)$, $t\in \mathbb R$, is not unpredictable.
\end{lemma}

\noindent \textbf{Proof.} Assume that the function $\phi$ is unpredictable. Then, there exists a sequence $\{\mu_n\}_{n\in\mathbb N}$ of real numbers, which diverges to infinity, such that $\left\| \phi(t+\mu_n) -\phi(t)\right\| \to 0$ as $n \to \infty$ uniformly on compact subsets of the real axis.

Firstly, suppose that $\left\| \xi(t_0)\right\| \geq 4H$ for some $t_0 \in\mathbb R$. One can confirm the existence of a natural number $n_1$ such that 
\begin{eqnarray} \label{discproofinq11}
\left\| \phi(t_0+ \mu_{n_1}) - \phi(t_0) \right\| <H
\end{eqnarray}
and
\begin{eqnarray} \label{discproofinq12}
\left\| \xi(t_0 + \mu_{n_1}) \right\| \leq H.
\end{eqnarray}
Utilizing (\ref{discproofinq12}) we attain that
\begin{eqnarray*}
\left\| \phi(t_0 + \mu_{n_1}) - \phi(t_0) \right\| \geq \left\| \xi(t_0) \right\| - \left\| \psi(t_0 + \mu_{n_1}) \right\|  - \left\| \psi(t_0) \right\| -  \left\| \xi(t_0 + \mu_{n_1}) \right\| \geq H.
\end{eqnarray*}
The last inequality contradicts to (\ref{discproofinq11}).

Secondly, let us consider the case $\left\| \xi(\theta_{k_0}+)\right\| > 4H$ for some $k_0 \in \mathbb Z$. There exists a natural number $n_2$ such that $\left\| \phi(t + \mu_{n_2}) - \phi(t) \right\| <H$ for $t \in [\theta_{k_0}, \theta_{k_0}+1]$ and $\left\|\xi(t+\mu_{n_2}) \right\| <H$ for $t\geq \theta_{k_0}.$ Accordingly, we have 
\begin{eqnarray*} \label{discproofinq13}
\left\| \phi((\theta_{k_0} + \mu_{n_2})+) - \phi(\theta_{k_0} +) \right\| \leq H
\end{eqnarray*} 
and 
\begin{eqnarray} \label{discproofinq14}
\left\| \xi((\theta_{k_0} + \mu_{n_2})+)\right\| \leq H.
\end{eqnarray} 
It can be obtained by means of (\ref{discproofinq14}) that
\begin{eqnarray*}
\left\| \phi((\theta_{k_0} + \mu_{n_2})+) - \phi(\theta_{k_0} +) \right\| \geq \left\| \xi(\theta_{k_0} +) \right\| - \left\| \psi((\theta_{k_0} + \mu_{n_2})+) \right\|  - \left\|\psi(\theta_{k_0} +) \right\| - \left\| \xi((\theta_{k_0} + \mu_{n_2})+)\right\| > H.
\end{eqnarray*}
This is also a contradiction. Thus, $\phi$ is not unpredictable. $\square$

%%%%%%%%%%%%

\section{Asymptotic Unpredictability in Impulsive Systems} \label{sec3}

Our purpose in this section is to investigate asymptotically unpredictable solutions of quasilinear impulsive systems of the form
\begin{eqnarray} \label{mainimpusyst}
&& x'(t) = Ax(t) + f (t,x(t)) + g_1(t)+g_2(t), \ t \neq \theta_k, \nonumber \\
&& \Delta x \big{|}_{t=\theta_k} = Bx(\theta_k) + h (x(\theta_k)) + \gamma_k,
\end{eqnarray}
where the $m\times m$ constant matrices $A$ and $B$ commute,  the functions $f:\mathbb R \times \mathbb R^m \to \mathbb R^m$, $h:\mathbb R^m \to \mathbb R^m$ are continuous in all their arguments, the sequence $\theta=\{\theta_k\}_{k\in\mathbb Z}$ of impulse moments belongs to the set $\mathcal{S}(p,\omega)$ for some fixed numbers $p\in\mathbb N$ and $\omega >0$, $\{\gamma_k\}_{k\in\mathbb Z}$ is a bounded sequence with $\displaystyle \lim_{k \to \infty} \left\|\gamma_k \right\| =0$, and $\Delta x \big{|}_{t=\theta_k} = x(\theta_k+) - x(\theta_k)$ for each $k\in\mathbb Z$ in which $x(\theta_k+) = \displaystyle \lim_{t \to \theta_k^+} x(t)$. Moreover, $g_1 \in \mathcal{BPC}(\theta)$ is the function satisfying the equation
\begin{eqnarray} \label{defnofg}
g_1(t) = \sigma_k,  
\end{eqnarray}
for $\theta_{kp}<t\leq \theta_{(k+1)p},$ $k\in\mathbb Z$, where $\{\sigma_k\}_{k\in\mathbb Z}$ is a fixed bounded sequence in $\mathbb R^m$, and $g_2 \in \mathcal{BPC}(\theta)$ is a function with $\displaystyle \lim_{t \to \infty} \left\|g_2(t) \right\| =0$.

In our research, we need the following assumptions on the impulsive system (\ref{mainimpusyst}).
\begin{itemize}
\item[\textbf{(A1)}] $\det(I+B)\neq0$, where $I$ denotes the $m \times m$ identity matrix,
\item[\textbf{(A2)}] The real parts of all eigenvalues of the matrix $A + \displaystyle \frac{p}{\omega} \textrm{Log}(I+B)$ are negative,
\item[\textbf{(A3)}] $f(t+\omega,x) = f(t,x)$ for each $(t,x) \in \mathbb R \times \mathbb R^m$,
\item[\textbf{(A4)}] The inequalities
$$\displaystyle \sup_{(t,x) \in [0,\omega) \times \mathbb R^m} \left\|f(t,x) \right\| \leq M_f$$
and
$$\displaystyle \sup_{x \in \mathbb R^m} \left\|h(x) \right\| \leq M_h$$
hold for some positive numbers $M_f$ and $M_h$,
\item[\textbf{(A5)}] The inequalities
$$ \left\|f(t,x_1) - f(t,x_2)\right\| \leq L_f \left\|x_1-x_2 \right\|, \ t\in\mathbb R, x_1,x_2 \in\mathbb R^m $$
and
$$\left\|h(x_1)-h(x_2) \right\| \leq L_h  \left\|x_1-x_2 \right\|, \  x_1,x_2 \in\mathbb R^m $$
hold for some positive numbers $L_f$ and $L_h$.
\end{itemize}

In what follows, for a given interval $J$ of the real axis, $i(J)$ stands for  the number of the terms of the sequence $\theta=\left\lbrace \theta_k\right\rbrace_{k\in\mathbb Z}$ that belong to $J$. For any real numbers $t_1$ and $t_2$ with $t_1<t_2$, the inequality
\begin{eqnarray} \label{ineqii}
i((t_1,t_2)) \leq p + \displaystyle \frac{p}{\omega} (t_2-t_1)
\end{eqnarray} 
is valid \cite{Fen2023}. Besides, one can attain under the assumption $(A1)$ that 
$$U(t,s) = e^{A(t-s)} (I+B)^{i([s,t))}, \ t>s$$
and $U(s,s)=I$, where $U(t,s)$ denotes the matriciant of the linear homogeneous impulsive system
\begin{eqnarray} \label{linimpusyst}
	&& x'(t) = Ax(t), \ t \neq \theta_k, \nonumber \\
	&& \Delta x \big{|}_{t=\theta_k} = Bx(\theta_k).
\end{eqnarray}
If the assumption $(A2)$ additionally holds, then according to the results of the books \cite{Akhmet2010,Samoilenko1995}, there exist numbers $N \geq 1$ and $\lambda>0$ such that 
\begin{eqnarray} \label{ineqmatriciant}
\left\| U(t,s) \right\| \leq N e^{-\lambda (t-s)}, \ t \geq s.
\end{eqnarray}  
The following assumptions are also required to approve asymptotic unpredictability.
\begin{itemize}
\item[\textbf{(A6)}] $N\left( \displaystyle \frac{L_f}{\lambda} +  \frac{p L_h}{1-e^{-\lambda \omega}}\right)<1$,
\item[\textbf{(A7)}] $NL_f + \displaystyle \frac{p}{\omega} \ln (1+NL_h)<\lambda$, 
\item[\textbf{(A8)}] $L_h\left\|(I+B)^{-1} \right\| <1$.
\end{itemize}

We utilize the following definitions of unpredictable and asymptotically unpredictable sequences \cite{Fen2018,Fen2024b}.

\begin{definition} (\cite{Fen2018}) \label{unpseq1}
A bounded sequence $\{\alpha_k\}_{k\in\mathbb Z}$ in $\mathbb R^m$ is called unpredictable if there exist a positive number $\upsilon_0$ (the unpredictability constant) and sequences $\{q_n\}_{n\in\mathbb N}$, $\{r_n\}_{n\in\mathbb N}$ of positive integers both of which diverge to infinity such that $\left\|\alpha_{k+q_n}-\alpha_k \right\| \to 0$ as $n \to\infty$ for each $k$  in bounded intervals of integers and $\left\| \alpha_{q_n+r_n} - \alpha_{r_n}\right\| \geq \upsilon_0$ for each $n\in\mathbb N$. 
\end{definition}

\begin{definition} (\cite{Fen2024b}) \label{asympunpseq1}
A bounded sequence $\{\sigma_k\}_{k\in\mathbb Z}$ in $\mathbb R^m$ is called asymptotically unpredictable if there exist an unpredictable sequence $\{\alpha_k\}_{k\in\mathbb Z}$ and a sequence $\{\beta_k\}_{k\in\mathbb Z}$ satisfying $\displaystyle \lim_{k \to \infty} \left\| \beta_k\right\|=0$ such that $\sigma_k = \alpha_k + \beta_k$ for every $k \in \mathbb Z$.
\end{definition}

The main result of the present study is mentioned in the subsequent theorem.
\begin{theorem} \label{newmainresult}
Suppose that the assumptions $(A1)-(A8)$ are valid. If the sequence $\{\sigma_k\}_{k\in\mathbb Z}$ is asymptotically unpredictable, then the impulsive system (\ref{mainimpusyst}) possesses a unique asymptotically unpredictable solution. Moreover, this solution is asymptotically stable.
\end{theorem}

\noindent \textbf{Proof.} Because the sequence $\{\sigma_k\}_{k \in \mathbb Z}$ is asymptotically unpredictable, there exist an unpredictable sequence $\{\alpha_k\}_{k\in\mathbb Z}$ in $\mathbb R^m$ and a sequence $\{\beta_k\}_{k\in\mathbb Z}$ in $\mathbb R^m$ satisfying $\left\| \beta_k\right\| \to 0$ as $k \to \infty$ such that $\sigma_k = \alpha_k + \beta_k$ for each $k\in\mathbb Z$.
 
Referring to Theorem $34$ and Theorem $89$ proposed in the book \cite{Samoilenko1995}, under the assumptions $(A1)-(A6)$, system (\ref{mainimpusyst}) has a unique bounded solution $\phi \in \mathcal{BPC}(\theta)$, which is asymptotically stable and satisfies the relation
\begin{eqnarray*}
\displaystyle \phi(t) = \int_{-\infty}^t U(t,s) \left[ f(s,\phi(s)) + g_1(s) + g_2(s)\right] ds + \sum_{-\infty < \theta_k <t} U(t,\theta_k+) \left[ h(\phi(\theta_k)) + \gamma_k\right] . 
\end{eqnarray*} 
To complete the proof, we will verify that $\phi(t)$ is asymptotically unpredictable. For that purpose, let us take into account the impulsive system
\begin{eqnarray} \label{impusyst2}
&& y'(t) = Ay(t) + f (t,y(t)) + g(t), \ t \neq \theta_k, \nonumber \\
&& \Delta y \big{|}_{t=\theta_k} = By(\theta_k) + h (y(\theta_k)),
\end{eqnarray}
where the piecewise constant function $g(t)$ is defined via the equation $g(t) = \alpha_k$ for  $\theta_{kp}<t\leq \theta_{(k+1)p},$ $k\in\mathbb Z$. The matrices $A$, $B$, the functions $f$, $h$, and the sequence $\theta=\{\theta_k\}_{k\in\mathbb Z}$ of impulse moments in (\ref{impusyst2}) are the same with the ones in (\ref{mainimpusyst}).  
According to Theorem 3.1 given in paper \cite{Fen2023}, if the assumptions $(A1)-(A8)$ are valid, then system 
(\ref{mainimpusyst}) possesses a unique unpredictable solution $\psi(t)$ which satisfies the equation
\begin{eqnarray*}
\displaystyle \psi(t) = \int_{-\infty}^t U(t,s) \left[ f(s,\psi(s)) + g(s)\right] ds + \sum_{-\infty < \theta_k <t} U(t,\theta_k+) h(\psi(\theta_k)). 
\end{eqnarray*}
To deduce the asymptotic unpredictability of $\phi(t)$, it is sufficient to show that $\left\|\phi(t) -\psi(t) \right\| \to 0$ as $t\to\infty$ in accordance with Definition \ref{impuldefn2}.

In the rest of the proof, we will take advantage of the notations 
$$d_0= \lambda-NL_f-\frac{p}{\omega}\ln(1+NL_h), \ M_{\beta} = \displaystyle \sup_{k \in \mathbb Z} \left\|\beta_k \right\|, \ M_{\gamma} = \displaystyle \sup_{k \in \mathbb Z} \left\|\gamma_k \right\|, \ M_{g_2} = \displaystyle \sup_{t \in \mathbb R} \left\|g_2(t) \right\|. $$
Let a positive number $\varepsilon$ be given, and fix a positive number $c_0$ satisfying $$c_0 \leq \displaystyle \frac{1}{H_1 + H_2},$$ where 
\begin{eqnarray} \label{numberh1}
H_1 = \left[  \frac{N}{\lambda} \left( 2M_f + M_{\beta} + M_{g_2} \right) + \frac{N p (2M_h + M_{\gamma})}{1-e^{-\lambda \omega}}  \right] (1+N L_h)^p 
\end{eqnarray}  
and
\begin{eqnarray} \label{numberh2}
H_2 = \left( \frac{2}{\lambda} + \frac{p}{1-e^{-\lambda \omega}} \right) \left[ N + \frac{N^2 L_f (1+NL_{h})^p}{d_0} + \frac{N^2 p L_{h} (1+NL_{h})^p}{1-e^{-d_0 \omega}} \right].  
\end{eqnarray}
Since both $\left\| \beta_k \right\|$ and $\left\| \gamma_k \right\|$ converge to zero as $k \to \infty$, there exist integers $k_1$ and $k_2$ such that 
\begin{eqnarray} \label{prooff1}
\left\| \beta_k \right\| < c_0 \varepsilon, \  k \geq k_1
\end{eqnarray} 
and 
\begin{eqnarray} \label{prooff2}
\left\| \gamma_k \right\| < c_0 \varepsilon, \  k \geq k_2.
\end{eqnarray}
With regard to (\ref{prooff1}) we have 
\begin{eqnarray} \label{prooff3}
\left\| g_1(t)-g(t)\right\| <c_0 \varepsilon
\end{eqnarray}
whenever $t > \theta_{k_{1} p}$. On the other hand, since $\left\| g_2(t)\right\| \to 0$ as $t\to\infty$, there is a real number $\widetilde{t}$ such that if $t \geq \widetilde{t}$, then $\left\| g_2(t)\right\| < c_0\varepsilon$. Let us denote
$$t_0 = \max \left\lbrace \theta_{k_{1} p}, \theta_{k_2}, \widetilde{t} \,  \right\rbrace.$$

Making use of the equation
\begin{eqnarray*}
\phi(t) - \psi(t) &=& \int_{-\infty}^t U(t,s) \left[f(s,\phi(s))- f(s,\psi(s)) + g_1(s) + g_2(s)-g(s) \right] ds  \\
&+&  \sum_{-\infty < \theta_k <t} U(t,\theta_k+) \left[ h(\phi(\theta_k)) -h(\psi(\theta_k)) + \gamma_k \right]  
\end{eqnarray*}
together with (\ref{prooff2}) and (\ref{prooff3}), we attain for $t >t_0$ that
\begin{eqnarray*}
\left\| \phi(t) - \psi(t) \right\| & \leq & \displaystyle \int_{-\infty}^{t_0} N(2M_f + M_{\beta} + M_{g_2}) e^{-\lambda (t-s)} ds 
+ \displaystyle \int_{t_0}^{t} NL_f e^{-\lambda (t-s)} \left\| \phi(s) - \psi(s) \right\| ds \\
&+&  \displaystyle \int_{t_0}^{t}  2 N c_0 \varepsilon e^{-\lambda (t-s)} ds + \sum_{-\infty < \theta_k \leq t_0} N (2M_h + M_{\gamma}) e^{-\lambda (t-\theta_k)} \\
&+& \displaystyle \sum_{t_0 < \theta_k < t} N L_{h} e^{-\lambda (t-\theta_k)} \left\| \phi(\theta_k) - \psi(\theta_k)\right\| +  \sum_{t_0 < \theta_k < t} N c_0 \varepsilon  e^{-\lambda (t-\theta_k)}.
\end{eqnarray*}
One can confirm that the estimates
$$\displaystyle \sum_{-\infty < \theta_k \leq t_0} e^{-\lambda (t-\theta_k)} \leq \frac{p}{1-e^{-\lambda \omega}} e^{-\lambda (t-t_0)} $$
and 
$$\displaystyle \sum_{t_0 < \theta_k < t} e^{-\lambda (t-\theta_k)} \leq \frac{p}{1-e^{-\lambda \omega}} \left( 1-e^{-\lambda (t-t_0+\omega)} \right)  $$
are satisfied. Therefore,
\begin{eqnarray} \label{ineqqproof1}
\left\| \phi(t) - \psi(t) \right\| & \leq & \frac{N}{\lambda} (2M_f + M_{\beta} + M_{g_2}) e^{-\lambda (t-t_0)} + \frac{2Nc_0\varepsilon}{\lambda} \left(1-e^{-\lambda (t-t_0)} \right) \nonumber \\
&+& \frac{N p (2M_h + M_{\gamma})}{1-e^{-\lambda \omega}} e^{-\lambda (t-t_0)} + \frac{N p c_{0} \varepsilon}{1-e^{-\lambda \omega}} \left( 1- e^{-\lambda (t-t_0 + \omega)}\right) \nonumber \\
&+& \displaystyle \int_{t_0}^{t} N L_{f} e^{-\lambda (t-s)}  \left\| \phi(s) - \psi(s) \right\| ds \nonumber \\
&+& \sum_{t_0 < \theta_k < t} N L_{h} e^{-\lambda (t-\theta_k)}  \left\| \phi(\theta_k) - \psi(\theta_k) \right\|
\end{eqnarray}
provided that $t >t_0$. 

Now, let us denote $$u(t) = e^{\lambda t} \left\| \phi(t) - \psi(t) \right\|$$ and 
$$\Gamma = \left[ \frac{N}{\lambda} (2 M_f + M_{\beta} + M_{g_2}) + \frac{Np (2 M_h + M_{\gamma})}{1-e^{-\lambda \omega}} - \frac{2Nc_0 \varepsilon}{\lambda} - \frac{N p c_{0} \varepsilon e^{-\lambda \omega}}{1-e^{-\lambda \omega}}  \right] e^{\lambda t_0}.$$
The inequality (\ref{ineqqproof1}) yields
\begin{eqnarray*}
u(t) \leq \Gamma + \left( \frac{2}{\lambda} + \frac{p}{1-e^{-\lambda \omega}} \right) N c_{0} \varepsilon e^{\lambda t} + \displaystyle \int_{t_{0}}^{t} N L_{f} u(s) ds +  \sum_{t_{0} < \theta_{k} < t} N L_{h} u(\theta_{k}), \ t>t_0. 
\end{eqnarray*}
Applying the Gronwall type inequality for piecewise continuous functions given in Theorem 16.2 of the book \cite{Bainov1992}, one can obtain that
\begin{eqnarray*}
u(t) & \leq & \Gamma + \left( \frac{2}{\lambda} + \frac{p}{1-e^{-\lambda \omega}} \right) N c_{0} \varepsilon e^{\lambda t} \\
&+& \Gamma N L_{f} \int_{t_0}^t  \left(1 + N L_h \right)^{i((s,t))} e^{NL_f(t-s)} ds \\
&+& \Gamma N L_h \sum_{t_0 < \theta_k < t}  \left(1+NL_h \right)^{i((\theta_k,t))} e^{NL_f(t-\theta_k)} \\
&+& \left(\frac{2}{\lambda} + \frac{p}{1-e^{-\lambda \omega}} \right) N^{2} L_f c_0 \varepsilon e^{N L_f t} \int_{t_0}^t \left(1+NL_h \right)^{i((s,t))}  e^{(\lambda-NL_f)s} ds \\
&+& \left(\frac{2}{\lambda} + \frac{p}{1-e^{-\lambda \omega}} \right) N^{2} L_h c_0 \varepsilon e^{N L_f t} \sum_{t_0 < \theta_k < t} \left(1+NL_h \right)^{i((\theta_k,t))}  e^{(\lambda-NL_f) \theta_k}.
\end{eqnarray*}
Using the estimate (\ref{ineqii}) and the identity
\begin{eqnarray*}
&& 1+\displaystyle \int_{t_0}^{t} N L_{f} (1+NL_h)^{i((s,t))} e^{NL_f(t-s)} ds + \sum_{t_0 < \theta_k < t} NL_{h} (1+NL_h)^{i((\theta_k,t))} e^{NL_f (t-\theta_k)} \\
&& = (1+NL_h)^{i((t_0,t))} e^{NL_f(t-t_0)},
\end{eqnarray*}
we attain that
\begin{eqnarray*}
u(t) & \leq  & \Gamma (1+NL_h)^p e^{(\lambda -d_0) (t-t_0)} + \left(\frac{2}{\lambda} + \frac{p}{1-e^{-\lambda \omega}} \right) N c_0 \varepsilon e^{\lambda t} \\
& + & \left(\frac{2}{\lambda} + \frac{p}{1-e^{-\lambda \omega}} \right) \frac{N^2 L_f (1+N L_h)^p c_0 \varepsilon}{d_0} e^{\lambda t} \left( 1- e^{-d_0 (t-t_0)}\right)  \\
& + & \left(\frac{2}{\lambda} + \frac{p}{1-e^{-\lambda \omega}} \right) \dfrac{N^2 p L_h (1+NL_h)^p c_0 \varepsilon }{1-e^{-d_0 \omega}} e^{\lambda t} \left(1-e^{-d_0 (t-t_0+\omega)} \right). 
\end{eqnarray*}
In accordance with the last inequality we have $$\left\|\phi(t) -\psi(t) \right\| <H_1e^{-d_0(t-t_0)} + H_2 c_0 \varepsilon$$ for $t>t_0$, in which $H_1$ and $H_2$ are respectively defined by (\ref{numberh1}) and (\ref{numberh2}). For that reason, if $t > t_0 + \Delta_0$, then $$\left\|\phi(t) -\psi(t) \right\| < (H_1 + H_2)c_0 \varepsilon \leq \varepsilon,$$ where $\Delta_0$ is a fixed positive number with $$\Delta_0 \geq \frac{1}{d_0} \ln \left( \frac{1}{c_0 \varepsilon}\right).$$ Hence, $\left\|\phi(t) -\psi(t) \right\| \to 0$ as $t \to \infty$. Consequently, there is a unique asymptotically unpredictable solution of the impulsive system (\ref{mainimpusyst}), which is asymptotically stable. $\square$

\begin{remark} The result of Theorem \ref{newmainresult} holds even if we replace the term $f(t,x(t))$ in (\ref{mainimpusyst})  with a function of the form $\widetilde{f}(x(t))$ and take into consideration the counterparts of the boundedness and Lipschitz conditions respectively mentioned in $(A4)$ and $(A5)$. \end{remark}

\section{Examples} \label{sec4}

This section is devoted to the exemplification of piecewise continuous asymptotically unpredictable functions and an impulsive system possessing such a solution. 

\subsection{Examples of Asymptotically Unpredictable Functions with Discontinuities} \label{examplesec41}

It was demonstrated in paper \cite{Fen2017} that the logistic map
\begin{eqnarray} \label{logisticmap}
z_{k+1} = 3.9 z_k (1-z_k), \ k\in\mathbb Z,
\end{eqnarray}
possesses an unpredictable solution $\{\alpha^*_k\}_{k\in\mathbb Z}$, which lies inside the unit interval $[0,1]$. In accordance with Definition \ref{unpseq1}, there exist a positive number $\upsilon_0$ and sequences $\{q_n\}_{n\in\mathbb N}$, $\{r_n\}_{n\in\mathbb N}$ of positive integers both of which diverge to infinity such that $|\alpha^*_{k+q_n}-\alpha^*_k| \to 0$ as $n \to\infty$ for each $k$ in bounded intervals of integers and $| \alpha^*_{q_n+r_n} - \alpha^*_{r_n}| \geq \upsilon_0$ for each $n\in\mathbb N$. 

Suppose that $\theta=\{\theta_k\}_{k\in\mathbb Z}$ is the sequence belonging to $\mathcal{S}\left(2, 3 \pi/2 \right)$ defined by 
\begin{eqnarray} \label{seqthetak}
\theta_k = \displaystyle \frac{\pi}{4} \left(3k + (-1)^k\right), \ k\in\mathbb Z.
\end{eqnarray}
Let us take into account the function $\xi_0:\mathbb R \to \mathbb R$ with $\xi_0(t) = \displaystyle \frac{13}{2k^2+1}$ for $t \in (\theta_{2k}, \theta_{2k+2}]$, $k\in\mathbb Z$. We define the piecewise continuous function $\phi:\mathbb R \to \mathbb R^3$ by means of the equation
\begin{eqnarray} \label{examplephi}
\phi(t) = \psi(t) + \xi(t),
\end{eqnarray}
where the functions $\psi:\mathbb R \to \mathbb R^3$ and $\xi:\mathbb R \to \mathbb R^3$ respectively satisfy the equations
\begin{eqnarray*}
\psi(t) = 
\begin{pmatrix} 
\displaystyle \frac{\alpha^*_k}{2} \\ -\alpha^*_k \\  \displaystyle \frac{\alpha^*_k}{4}
\end{pmatrix}
\end{eqnarray*}
for $t \in (\theta_{2k}, \theta_{2k+2}]$, $k\in\mathbb Z$, and
\begin{eqnarray*}
\xi(t) = 
\begin{pmatrix} 
\tanh t -1 \\ \xi_0 (t) \\  \displaystyle \frac{8}{3e^t+e^{-t}}
\end{pmatrix}.
\end{eqnarray*}

We will verify that $\phi$ is asymptotically unpredictable and that it is not unpredictable. Firstly, item $(i)$ of Definition \ref{impuldefn1} is fulfilled since $\psi$ is constant on each of the intervals $(\theta_{2k}, \theta_{2k+2}]$, $k\in\mathbb Z$. Now, take a positive number $\varepsilon$ and a compact subset $\mathcal{C}$ of $\mathbb R$. Fix integers $k_1$, $k_2$ with $k_2>k_1$ and $\mathcal{C} \subset \left( \theta_{2k_1}, \theta_{2k_2+2} \right]$. There is a natural number $n_0$ such that if $n\geq n_0$, then $\left|\alpha_{k+q_n} -\alpha_k \right|< 4\varepsilon/\sqrt{21}$ for $k=k_1,k_{1}+1,\ldots,k_2$. For each such value of $k$, if $t \in \left( \theta_{2k}, \theta_{2k+2} \right]$, then we have $$\left\| \psi(t + \mu_n) - \psi(t)\right\| = \displaystyle \frac{\sqrt{21}}{4} \left| \alpha_{k+q_n} - \alpha_k\right| <\varepsilon.$$ Accordingly, $\left\| \psi(t + \mu_n) - \psi(t)\right\| \to 0$ as $n \to \infty$ uniformly on $\mathcal{C}$. On the other hand, let us denote $\tau_n=\theta_{2r_n+1}$, $n \in \mathbb N$, and $\nu=\pi/8$. For each $t \in [\tau_n - \nu, \tau_n +\nu]$ and $n\in\mathbb N$, it can be shown that $$\left\| \psi(t + \mu_n) - \psi(t)\right\| > \frac{\sqrt{21}}{4} \upsilon_0.$$ This approves the unpredictability of $\psi$. Because $ \left\| \xi(t) \right\| \to 0$ as $t \to \infty$, the function $\phi$ is asymptotically unpredictable by Definition \ref{impuldefn2}. One can confirm utilizing Lemma \ref{disconlemma4} that $\phi$ is not unpredictable since $\displaystyle \sup_{t \in \mathbb R} \left\| \psi(t)\right\|  \leq \sqrt{21}/4$ and $\left\| \xi(0)\right\| > \sqrt{21}$. 

Finally, we define the functions 
$$\phi_1(t) = \Omega \phi(t) + \begin{pmatrix} 1 \\ -2 \\ -1 \end{pmatrix},$$ where $\Omega= \textrm{diag}(4,-1,2)$,
$$\phi_2(t) = \phi(t) + \begin{pmatrix} \arctan t \\ \displaystyle \frac{t^2}{5+t^2} \\ 1- \textrm{sech} (2t) \end{pmatrix}, $$ and $$\phi_3(t) = \phi(t+2).$$
The functions $\phi_1$, $\phi_2$, and $\phi_3$ are also asymptotically unpredictable by Lemma \ref{disconlemma1}, Lemma \ref{disconlemma2}, and Lemma \ref{disconlemma3}, respectively.

\subsection{An Impulsive System with an Asymptotically Unpredictable Solution}

Let us consider the impulsive system
\begin{eqnarray} \label{exampleimp}
&& x'_1 (t) = -4x_1(t) + x_2(t) + \frac{1}{9} \arctan(x_2(t)) + \Phi_1(t)+ \displaystyle \frac{8\cos t}{1+t^2},  \nonumber \\
&& x'_2 (t) = -2x_1(t) -6x_2(t) + \Phi_2(t) + \sin\left( \frac{4t}{3}\right) + 5\textrm{sech} t,  \  t \neq \theta_k,  \nonumber \\
&& \Delta x_1 \big{|}_{t=\theta_k} = - \frac{3}{5} x_1(\theta_k) +  \frac{1}{15}\sin(x_1(\theta_k)), \nonumber \\
&& \Delta x_2 \big{|}_{t=\theta_k} =  - \frac{3}{5} x_2(\theta_k) + \frac{1}{20} \cos(x_2(\theta_k))  + e^{-2k^2}, 
\end{eqnarray}
where the sequence $\theta=\{\theta_k\}_{k\in\mathbb Z}$ of impulse moments satisfies the equation (\ref{seqthetak}), $\Phi_1(t) = 1-\alpha^*_k -\tanh k$ and $\Phi_2(t) = \alpha^*_k + \displaystyle \frac{1}{1+e^k}$ for $t\in(\theta_{2k}, \theta_{2k+2}]$, $k\in\mathbb Z$, and $\{\alpha^*_k\}_{k\in\mathbb Z}$ is an unpredictable solution of (\ref{logisticmap}) that belongs to the interval $[0,1]$. System (\ref{exampleimp}) is in the form of (\ref{mainimpusyst})
with $$A=\begin{pmatrix} -4 && 1 \\ -2 && -6\end{pmatrix}, \  B=\begin{pmatrix} - \displaystyle \frac{3}{5} && 0 \\ 0 && - \displaystyle \frac{3}{5} \end{pmatrix},$$
$$f(t,x_1,x_2) = \begin{pmatrix} \displaystyle \frac{1}{9} \arctan(x_2) \\  \sin\left(\displaystyle \frac{4t}{3}\right) \end{pmatrix}, \ h(x_1,x_2) = \begin{pmatrix} \displaystyle \frac{1}{15} \sin(x_1) \vspace{.1cm} \\  \displaystyle \frac{1}{20} \cos(x_2) \end{pmatrix},$$
$$g_1(t) = \begin{pmatrix} \Phi_1(t) \\ \Phi_2(t)\end{pmatrix}, \ g_2(t) = \begin{pmatrix} \displaystyle \frac{8\cos t}{1+t^2}  \\  5\textrm{sech} t\end{pmatrix},$$
and $\gamma_k=e^{-2k^2}$, $k\in\mathbb Z$.

The matrix $A + \displaystyle \frac{p}{\omega} \textrm{Log}(I+B)$ admits the eigenvalues $-5+\displaystyle \frac{4}{3\pi} \ln \left( \displaystyle \frac{2}{5}\right)  \pm i$, where $p=2$ and $\omega=3\pi/2$. Additionally, the matriciant $U(t,s)$ of the linear homogeneous impulsive system 
\begin{eqnarray*} \label{exampleimp2}
	&& x'_1 (t) = -4x_1(t) + x_2(t),  \nonumber \\
	&& x'_2 (t) = -2x_1(t) -6x_2(t), \ t \neq \theta_k, \nonumber \\
	&& \Delta x_1 \big{|}_{t=\theta_k} = - \frac{3}{5} x_1(\theta_k),  \nonumber \\
	&& \Delta x_2 \big{|}_{t=\theta_k} =  - \frac{3}{5} x_2(\theta_k) 
\end{eqnarray*}
satisfies the equation 
$$
U(t,s) = \left( \frac{2}{5}\right)^{i([s,t))} e^{-5(t-s)} P \begin{pmatrix} \cos(t-s) && \sin(t-s) \\ -\sin(t-s) && \cos(t-s) \end{pmatrix} P^{-1}, \ t>s, 
$$ 
where $$P=\begin{pmatrix} 0  && 1 \\ 1 && -1 \end{pmatrix}.$$ 
Hence, (\ref{ineqmatriciant}) holds with $N=\left\| P\right\|  \left\| P^{-1}\right\| = \left( 3+\sqrt{5}\right) /2$ and $\lambda=5$. Taking $M_f=(\pi^2 +324)^{1/2}/18$, $M_h=1/12$, $L_f=1/9$, and $L_h=1/15$, it can be verified that the assumptions $(A1)-(A8)$ are valid for system (\ref{exampleimp}).  Besides, the equation $g_1(t) = \sigma_k$ is fulfilled for each $t\in(\theta_{2k}, \theta_{2k+2}]$ and $k\in\mathbb Z$, where $\{\sigma_k\}_{k\in\mathbb Z}$ is the sequence in $\mathbb R^2$ satisfying $\sigma_k=\alpha_k+\beta_k$, $k\in\mathbb Z$, in which 
\begin{eqnarray} \label{alphadefn1}
	\alpha_k= \begin{pmatrix} -\alpha^*_k \\ \alpha^*_k \end{pmatrix},
\end{eqnarray}
\begin{eqnarray} \label{betadefn1}
	\beta_k= \begin{pmatrix} 1-\tanh k \\ \displaystyle \frac{1}{1+e^k}\end{pmatrix}. 
\end{eqnarray}
The sequence $\{\sigma_k\}_{k\in\mathbb Z}$ is asymptotically unpredictable since $\{\alpha_k\}_{k\in\mathbb Z}$ is unpredictable and $\left\|\beta_k \right\|  \to 0$ as $k \to \infty$. Therefore, according to Theorem \ref{newmainresult}, there is a unique asymptotically unpredictable solution of the impulsive system (\ref{exampleimp}), which is asymptotically stable.
 
\begin{remark}
Let $\widetilde{\psi}:\mathbb R \to \mathbb R^2$ and $\overline{\xi}:\mathbb R \to \mathbb R^2$ be the functions  defined by $\widetilde{\psi}(t)=\alpha_k$ and $\overline{\xi}(t) = \beta_k$ for $(\theta_{2k}, \theta_{2k+2}]$, $k\in\mathbb Z$, where the sequences $\{\alpha_k\}_{k\in\mathbb Z}$ and $\{\beta_k\}_{k\in\mathbb Z}$ respectively satisfy (\ref{alphadefn1}) and (\ref{betadefn1}). Taking into consideration the functions $g_1$ and $g_2$ aforementioned in this section, one can confirm that $g_1(t)+g_2(t)=\widetilde{\psi}(t) + \widetilde{\xi}(t)$, in which $\widetilde{\xi}(t)=\overline{\xi}(t)+g_2(t)$. It can be shown using the arguments of Section \ref{examplesec41} that $\widetilde{\psi}$ is unpredictable. Moreover, $\big\| \widetilde{\xi}(t) \big\| \to 0$ as $t \to \infty$. For that reason $g_1(t) + g_2(t)$ is a piecewise continuous asymptotically unpredictable function. On the other hand, $\displaystyle \sup_{t\in\mathbb R} \big\|\widetilde{\psi}(t) \big\| \leq \sqrt 2$ and $\big\|\widetilde{\xi}(0)\big\| >4\sqrt 2$. Hence, $g_1(t) + g_2(t)$ is not unpredictable by Lemma \ref{disconlemma4}, and accordingly, system (\ref{exampleimp}) comprises a term that is asymptotically unpredictable but not unpredictable. This manifests the novelty of Theorem \ref{newmainresult} such that the existence and uniqueness of asymptotically unpredictable solutions of impulsive systems cannot be reduced to the case of unpredictable ones, which were discussed in papers \cite{Fen2023,Akhmet2022}. 
\end{remark}

\section{Conclusion} \label{concsec}

The present study is essentially concerned with asymptotically unpredictable motions comprising discontinuities. This kind of functions are newly introduced and related fundamental properties are obtained. It is rigorously proved in Lemma \ref{disconlemma4} that the set of piecewise continuous unpredictable functions is properly included in the set of asymptotically unpredictable ones. Asymptotically unpredictable motions generated by impulsive systems are also under discussion. A contribution to the qualitative theory of impulsive differential equations is performed such that the existence and uniqueness of asymptotically unpredictable solutions are demonstrated. Owing to the existence of piecewise continuous asymptotically unpredictable functions which are not unpredictable, the result of Theorem \ref{newmainresult} cannot be obtained from the study \cite{Fen2023}. This manifests the novelty of the theorem. On the other hand, system (\ref{mainimpusyst}) is different from the model discussed in \cite{Fen2023} since the term $g_1(t)$ is constructed through an asymptotically unpredictable sequence instead of an unpredictable one. In addition to that the differential equation and the impulse action in model (\ref{mainimpusyst}) are constructed respectively utilizing the function $g_2(t)$ and the sequence $\{\gamma_k \}_{k\in\mathbb Z}$, both of which converge to zero. It is worth noting that such terms do not take place in the model investigated in study \cite{Fen2023}. In the future, the results mentioned in this paper can be applied to epidemic dynamics \cite{Wang2024} and developed for impulsive fractional differential and evolution equations \cite{Wang2014,Shi2016}.

%\section*{CRediT authorship contribution statement}
%
%\textbf{Mehmet Onur Fen:} Conceptualization, Methodology, Validation, Formal analysis, Investigation, Resources, Writing-original draft, Writing-review \& editing.
%\textbf{Fatma Tokmak Fen:} Conceptualization, Methodology, Validation, Formal analysis, Investigation, Resources, Writing-original draft, Writing-review \& editing.

\section*{Declaration of competing interest}

The authors declare that they have no known competing financial interests or personal relationships that could have appeared to influence the work reported in this paper.

\section*{Data availability}

No data was used for the research described in the article.

%%%%%%%%%%%%%%%%%%% References %%%%%%%%%%%%%%%%%%%%%%%%%%%%%

\end{document}